\documentclass[reqno, 12pt]{amsart}
\pdfoutput=1
\makeatletter
\let\origsection=\section \def\section{\@ifstar{\origsection*}{\mysection}}
\def\mysection{\@startsection{section}{1}\z@{.7\linespacing\@plus\linespacing}{.5\linespacing}{\normalfont\scshape\centering\S}}
\makeatother

\usepackage{amsmath,amssymb,amsthm}
\usepackage{mathrsfs}
\usepackage{mathabx}\changenotsign
\usepackage{dsfont}
\usepackage{bbm}

\usepackage{xcolor}
\usepackage[backref=section]{hyperref}
\usepackage[ocgcolorlinks]{ocgx2}
\hypersetup{
	colorlinks=true,
	linkcolor={red!60!black},
	citecolor={green!60!black},
	urlcolor={blue!60!black},
}

\usepackage[open,openlevel=2,atend]{bookmark}

\usepackage[abbrev,msc-links,backrefs]{amsrefs}
\usepackage{doi}

\renewcommand{\PrintDOI}[1]{\doi{#1}}

\usepackage[T1]{fontenc}
\usepackage{lmodern}
\usepackage[babel]{microtype}
\usepackage[english]{babel}

\usepackage{geometry}
\numberwithin{equation}{section}
\numberwithin{figure}{section}

\usepackage{enumitem}

\let\polishlcross=\l
\def\l{\ifmmode\ell\else\polishlcross\fi}

\let\setminus=\smallsetminus

\makeatletter
\def\moverlay{\mathpalette\mov@rlay}
\def\mov@rlay#1#2{\leavevmode\vtop{   \baselineskip\z@skip \lineskiplimit-\maxdimen
		\ialign{\hfil$\m@th#1##$\hfil\cr#2\crcr}}}
\newcommand{\charfusion}[3][\mathord]{
	#1{\ifx#1\mathop\vphantom{#2}\fi
		\mathpalette\mov@rlay{#2\cr#3}
	}
	\ifx#1\mathop\expandafter\displaylimits\fi}
\makeatother

\newcommand{\dcup}{\charfusion[\mathbin]{\cup}{\cdot}}

\DeclareFontFamily{U}  {MnSymbolC}{}
\DeclareSymbolFont{MnSyC}         {U}  {MnSymbolC}{m}{n}
\DeclareFontShape{U}{MnSymbolC}{m}{n}{
	<-6>  MnSymbolC5
	<6-7>  MnSymbolC6
	<7-8>  MnSymbolC7
	<8-9>  MnSymbolC8
	<9-10> MnSymbolC9
	<10-12> MnSymbolC10
	<12->   MnSymbolC12}{}
\DeclareMathSymbol{\powerset}{\mathord}{MnSyC}{180}

\usepackage{tikz}
\usetikzlibrary{calc,decorations.pathmorphing}
\usetikzlibrary{arrows,decorations.pathreplacing}
\pgfdeclarelayer{background}
\pgfdeclarelayer{foreground}
\pgfdeclarelayer{front}
\pgfsetlayers{background,main,foreground,front}

\usepackage{multicol}
\usepackage{subcaption}
\newcommand{\pedge}[9]{
	
	\ifx\relax#6\relax
	\def\qoffs{0pt}
	\else
	\def\qoffs{#6}
	\fi
	
	\def\phedge{
		($#1+#5!\qoffs!-90:#2-#5$) -- 
		($#2+#1!\qoffs!-90:#3-#1$) -- 
		($#3+#2!\qoffs!-90:#4-#2$) -- 
		($#4+#3!\qoffs!-90:#5-#3$) -- 
		($#5+#4!\qoffs!-90:#1-#4$) -- cycle}

	\coordinate (12) at ($#1!\qoffs!90:#2$);
	\coordinate (15) at ($#1!\qoffs!-90:#5$);
	\coordinate (23) at ($#2!\qoffs!90:#3$);
	\coordinate (21) at ($#2!\qoffs!-90:#1$);
	\coordinate (34) at ($#3!\qoffs!90:#4$);
	\coordinate (32) at ($#3!\qoffs!-90:#2$);
	\coordinate (45) at ($#4!\qoffs!90:#5$);
	\coordinate (43) at ($#4!\qoffs!-90:#3$);
	\coordinate (51) at ($#5!\qoffs!90:#1$);
	\coordinate (54) at ($#5!\qoffs!-90:#4$);

	\def\nphedge{
		(15) let \p1=($(15)-#1$), \p2=($(12)-#1$) in 
		arc[start angle={atan2(\y1,\x1)}, delta angle={atan2(\y2,\x2)-atan2(\y1,\x1)-360*(atan2(\y2,\x2)-atan2(\y1,\x1)>0)}, x radius=\qoffs, y radius=\qoffs] --
		(21) let \p1=($(21)-#2$), \p2=($(23)-#2$) in 
		arc[start angle={atan2(\y1,\x1)}, delta angle={atan2(\y2,\x2)-atan2(\y1,\x1)-360*(atan2(\y2,\x2)-atan2(\y1,\x1)>0)}, x radius=\qoffs, y radius=\qoffs] --
		(32) let \p1=($(32)-#3$), \p2=($(34)-#3$) in 
		arc[start angle={atan2(\y1,\x1)}, delta angle={atan2(\y2,\x2)-atan2(\y1,\x1)-360*(atan2(\y2,\x2)-atan2(\y1,\x1)>0)}, x radius=\qoffs, y radius=\qoffs] --
		(43) let \p1=($(43)-#4$), \p2=($(45)-#4$) in 
		arc[start angle={atan2(\y1,\x1)}, delta angle={atan2(\y2,\x2)-atan2(\y1,\x1)-360*(atan2(\y2,\x2)-atan2(\y1,\x1)>0)}, x radius=\qoffs, y radius=\qoffs] --
		(54) let \p1=($(54)-#5$), \p2=($(51)-#5$) in 
		arc[start angle={atan2(\y1,\x1)}, delta angle={atan2(\y2,\x2)-atan2(\y1,\x1)-360*(atan2(\y2,\x2)-atan2(\y1,\x1)>0)}, x radius=\qoffs, y radius=\qoffs] --
		cycle}

	\ifx\relax#7\relax
	\def\plwidth{1pt}
	\else
	\def\plwidth{#7}
	\fi
	
	\ifx\relax#9\relax
	\fill \nphedge;
	\else
	\fill[#9]\nphedge;
	\fi
	
	\ifx\relax#8\relax
	\draw[line width=\plwidth,rounded corners=\qoffs]\nphedge;
	\else
	\draw[line width=\plwidth,#8]\nphedge;
	\fi
}

\newcommand{\qedge}[7]{
	
	\ifx\relax#4\relax
	\def\qoffs{0pt}
	\else
	\def\qoffs{#4}
	\fi
	
	\def\qhedge{
		($#1+#3!\qoffs!-90:#2-#3$) --
		($#2+#1!\qoffs!-90:#3-#1$) --
		($#3+#2!\qoffs!-90:#1-#2$) -- cycle}

	\coordinate (12) at ($#1!\qoffs!90:#2$);
	\coordinate (13) at ($#1!\qoffs!-90:#3$);
	\coordinate (23) at ($#2!\qoffs!90:#3$);
	\coordinate (21) at ($#2!\qoffs!-90:#1$);
	\coordinate (31) at ($#3!\qoffs!90:#1$);
	\coordinate (32) at ($#3!\qoffs!-90:#2$);
	
	\def\nqhedge{
		(13) let \p1=($(13)-#1$), \p2=($(12)-#1$) in
		arc[start angle={atan2(\y1,\x1)}, delta angle={atan2(\y2,\x2)-atan2(\y1,\x1)-360*(atan2(\y2,\x2)-atan2(\y1,\x1)>0)}, x radius=\qoffs, y radius=\qoffs] --
		(21) let \p1=($(21)-#2$), \p2=($(23)-#2$) in
		arc[start angle={atan2(\y1,\x1)}, delta angle={atan2(\y2,\x2)-atan2(\y1,\x1)-360*(atan2(\y2,\x2)-atan2(\y1,\x1)>0)}, x radius=\qoffs, y radius=\qoffs] --
		(32) let \p1=($(32)-#3$), \p2=($(31)-#3$) in
		arc[start angle={atan2(\y1,\x1)}, delta angle={atan2(\y2,\x2)-atan2(\y1,\x1)-360*(atan2(\y2,\x2)-atan2(\y1,\x1)>0)}, x radius=\qoffs, y radius=\qoffs] --
		cycle}
	
	\ifx\relax#5\relax
	\def\qlwidth{1pt}
	\else
	\def\qlwidth{#5}
	\fi
	
	\ifx\relax#7\relax
	\fill \nqhedge;
	\else
	\fill[#7]\nqhedge;
	\fi
	
	\ifx\relax#6\relax
	\draw[line width=\qlwidth,rounded corners=\qoffs]\nqhedge;
	\else
	\draw[line width=\qlwidth,#6]\nqhedge;
	\fi
}

\newcommand{\redge}[8]{
	
	\ifx\relax#5\relax
	\def\qoffs{0pt}
	\else
	\def\qoffs{#5}
	\fi
	
	\def\rhedge{
		($#1+#4!\qoffs!-90:#2-#4$) -- 
		($#2+#1!\qoffs!-90:#3-#1$) -- 
		($#3+#2!\qoffs!-90:#4-#2$) -- 
		($#4+#3!\qoffs!-90:#1-#3$) -- cycle}

	\coordinate (12) at ($#1!\qoffs!90:#2$);
	\coordinate (14) at ($#1!\qoffs!-90:#4$);
	\coordinate (23) at ($#2!\qoffs!90:#3$);
	\coordinate (21) at ($#2!\qoffs!-90:#1$);
	\coordinate (34) at ($#3!\qoffs!90:#4$);
	\coordinate (32) at ($#3!\qoffs!-90:#2$);
	\coordinate (41) at ($#4!\qoffs!90:#1$);
	\coordinate (43) at ($#4!\qoffs!-90:#3$);
	
	\def\nrhedge{
		(14) let \p1=($(14)-#1$), \p2=($(12)-#1$) in 
		arc[start angle={atan2(\y1,\x1)}, delta angle={atan2(\y2,\x2)-atan2(\y1,\x1)-360*(atan2(\y2,\x2)-atan2(\y1,\x1)>0)}, x radius=\qoffs, y radius=\qoffs] --
		(21) let \p1=($(21)-#2$), \p2=($(23)-#2$) in 
		arc[start angle={atan2(\y1,\x1)}, delta angle={atan2(\y2,\x2)-atan2(\y1,\x1)-360*(atan2(\y2,\x2)-atan2(\y1,\x1)>0)}, x radius=\qoffs, y radius=\qoffs] --
		(32) let \p1=($(32)-#3$), \p2=($(34)-#3$) in 
		arc[start angle={atan2(\y1,\x1)}, delta angle={atan2(\y2,\x2)-atan2(\y1,\x1)-360*(atan2(\y2,\x2)-atan2(\y1,\x1)>0)}, x radius=\qoffs, y radius=\qoffs] --
		(43) let \p1=($(43)-#4$), \p2=($(41)-#4$) in 
		arc[start angle={atan2(\y1,\x1)}, delta angle={atan2(\y2,\x2)-atan2(\y1,\x1)-360*(atan2(\y2,\x2)-atan2(\y1,\x1)>0)}, x radius=\qoffs, y radius=\qoffs] --
		cycle}
	
	\ifx\relax#6\relax
	\def\rlwidth{1pt}
	\else
	\def\rlwidth{#6}
	\fi
	
	\ifx\relax#8\relax
	\fill \nrhedge;
	\else
	\fill[#8]\nrhedge;
	\fi
	
	\ifx\relax#7\relax
	\draw[line width=\rlwidth,rounded corners=\qoffs]\nrhedge;
	\else
	\draw[line width=\rlwidth,#7]\nrhedge;
	\fi
}

\let\epsilon=\varepsilon

\let\rho=\varrho
\let\theta=\vartheta

\newtheoremstyle{note}  {4pt}  {4pt}  {\sl}  {}  {\bfseries}  {.}  {.5em}          {}
\newtheoremstyle{introthms}  {3pt}  {3pt}  {\itshape}  {}  {\bfseries}  {.}  {.5em}          {\thmnote{#3}}
\newtheoremstyle{remark}  {2pt}  {2pt}  {\rm}  {}  {\bfseries}  {.}  {.3em}          {}

\theoremstyle{plain}
\newtheorem{theorem}{Theorem}[section]

\theoremstyle{note}

\theoremstyle{remark}

\usepackage{lineno}
\newcommand*\patchAmsMathEnvironmentForLineno[1]{
	\expandafter\let\csname old#1\expandafter\endcsname\csname #1\endcsname
	\expandafter\let\csname oldend#1\expandafter\endcsname\csname end#1\endcsname
	\renewenvironment{#1}
	{\linenomath\csname old#1\endcsname}
	{\csname oldend#1\endcsname\endlinenomath}}
\newcommand*\patchBothAmsMathEnvironmentsForLineno[1]{
	\patchAmsMathEnvironmentForLineno{#1}
	\patchAmsMathEnvironmentForLineno{#1*}}
\AtBeginDocument{
	\patchBothAmsMathEnvironmentsForLineno{equation}
	\patchBothAmsMathEnvironmentsForLineno{align}
	\patchBothAmsMathEnvironmentsForLineno{flalign}
	\patchBothAmsMathEnvironmentsForLineno{alignat}
	\patchBothAmsMathEnvironmentsForLineno{gather}
	\patchBothAmsMathEnvironmentsForLineno{multline}
}

\def\ex{\text{\rm ex}}

\usepackage{scalerel}

\makeatletter
\newcommand{\overrighharpoonup}[1]{\ThisStyle{%
		\vbox {\m@th\ialign{##\crcr
				\rightharpoonupfill \crcr
				\noalign{\kern-\p@\nointerlineskip}
				$\hfil\SavedStyle#1\hfil$\crcr}}}}

\def\rightharpoonupfill{%
	$\SavedStyle\m@th\mkern+0.8mu\cleaders\hbox{$\shortbar\mkern-4mu$}\hfill\rightharpoonuptip\mkern+0.8mu$}

\def\rightharpoonuptip{%
	\raisebox{\z@}[2pt][1pt]{\scalebox{0.55}{$\SavedStyle\rightharpoonup$}}}

\def\shortbar{%
	\smash{\scalebox{0.55}{$\SavedStyle\relbar$}}}
\makeatother

\makeatletter
\newcommand{\overlefharpoonup}[1]{\ThisStyle{%
		\vbox {\m@th\ialign{##\crcr
				\leftharpoonupfill \crcr
				\noalign{\kern-\p@\nointerlineskip}
				$\hfil\SavedStyle#1\hfil$\crcr}}}}

\def\leftharpoonupfill{%
	$\SavedStyle\m@th\mkern+0.8mu\cleaders\hbox{$\shortbar\mkern-4mu$}\hfill\leftharpoonuptip\mkern+0.8mu$}

\def\leftharpoonuptip{%
	\raisebox{\z@}[2pt][1pt]{\scalebox{0.55}{$\SavedStyle\leftharpoonup$}}}

\makeatletter
\newsavebox\myboxA
\newsavebox\myboxB
\newlength\mylenA

\newcommand*\xoverline[2][0.75]{%
	\sbox{\myboxA}{$\m@th#2$}%
	\setbox\myboxB\null
	\ht\myboxB=\ht\myboxA%
	\dp\myboxB=\dp\myboxA%
	\wd\myboxB=#1\wd\myboxA
	\sbox\myboxB{$\m@th\overline{\copy\myboxB}$}
	\setlength\mylenA{\the\wd\myboxA}
	\addtolength\mylenA{-\the\wd\myboxB}%
	\ifdim\wd\myboxB<\wd\myboxA%
	\rlap{\hskip 0.5\mylenA\usebox\myboxB}{\usebox\myboxA}%
	\else
	\hskip -0.5\mylenA\rlap{\usebox\myboxA}{\hskip 0.5\mylenA\usebox\myboxB}%
	\fi}
\makeatother

\begin{document}
	
	\title[Hypergraphs accumulate]
	{Hypergraphs accumulate}
	
	\author[D.~Conlon]{David Conlon}
	\address{Department of Mathematics, California Institute of Technology, USA}
	\email{\{dconlon, schuelke\}@caltech.edu}
	
	\author[B.~Sch\"ulke]{Bjarne Sch\"ulke}

	
	\begin{abstract}
        We show that for every integer~$k\geq3$, the set of Tur\'an densities of~$k$-uniform hypergraphs 
        has an accumulation point in~$[0,1)$. 
        In particular, $1/2$ is an accumulation point for the set of Tur\'an densities of~$3$-uniform hypergraphs.
	\end{abstract}
	
	\maketitle
	
	\section{Introduction}
        
        For~$k\in\mathds{N}$, a~$k$-uniform hypergraph (or~$k$-graph)~$H=(V,E)$ consists of a vertex set~$V$ and an edge set~$E\subseteq V^{(k)}=\{e\subseteq V:\vert e\vert=k\}$.
        Given~$n\in\mathds{N}$ and a~$k$-graph~$F$, the extremal number~$\ex(n,F)$ is the maximum number of edges in a~$k$-graph~$H$ with~$n$ vertices that does not contain a copy of~$F$.
        The Tur\'an density of~$F$ is then given by
        $$\pi(F)=\lim_{n\to\infty}\frac{\ex(n,F)}{\binom{n}{k}},$$
        where the limit is known, by a simple monotonicity argument~\cite{KNS:64}, to be well-defined. The problem of determining these Tur\'an densities is one of the oldest and most fundamental questions in extremal combinatorics.
        
        When~$k=2$, that is, if~$F$ is a graph,~$\pi(F)$ is essentially completely understood, with the final result, the culmination of work by Tur\'an~\cite{T:41}, Erd\H{o}s and Stone~\cite{ES:46} and Erd\H{o}s and Simonovits~\cite{ES:66}, 
        saying that~$\pi(F)=\frac{\chi(F)-2}{\chi(F)-1}$, where~$\chi(F)$ is the chromatic number of~$F$.
        In contrast, very little is known about Tur\'an densities for~$k\geq 3$, 
        with even the problem of determining the Tur\'an density of the complete~$3$-graph on four vertices, a question first raised by Tur\'an in 1941~\cite{T:41}, remaining wide open. 
        For more on what is known about hypergraph Tur\'an densities, we refer the interested reader to the many surveys on the topic~\cites{F:91,K:11,S:95}.
        
        Despite the difficulty of determining the Tur\'an density of specific~$3$-graphs, one might try to study the distribution of the set~$\Pi^{(k)}=\{\pi(F):F\text{ is a }k\text{-graph}\}$ of Tur\'an densities of $k$-uniform hypergraphs. 
        For example, by a result of Erd\H{o}s~\cite{E:64} saying that~$\pi(F)=0$ if and only if~$F$ is a~$k$-partite~$k$-graph, we know that there is no~$k$-graph~$F$ with~$\pi(F)\in(0,k!/k^k)$. 
        However, this direction turned out to be similarly difficult and, beyond that simple result and the identification of some specific points in the set, very little is known about $\Pi^{(k)}$.
        
        If instead one considers 
        the set~$\Pi_{\infty}^{(k)}=\{\pi(\mathcal{F}):\mathcal{F}\text{ is a family of }k\text{-graphs}\}$, more is known. Of particular note here is the result of Frankl and R\"odl~\cite{FR:84} showing that~$\Pi^{(k)}_{\infty}$ is not well-ordered, thereby disproving the jumping conjecture, for which Erd\H{o}s had offered \$1000, saying that there is a non-trivial gap or jump between every two elements of~$\Pi^{(k)}_{\infty}$. 
        For more on the existence and non-existence of jumps, see, for instance,~\cite{BT:11,FPRT:07}.

        A more systematic study of~$\Pi^{(k)}_{\infty}$ and~$\Pi^{(k)}_{\text{fin}}=\{\pi(\mathcal{F}):\mathcal{F}\text{ is a finite family of }k\text{-graphs}\}$ was undertaken by Pikhurko~\cite{P:14}, who, roughly speaking, showed that for every iterative blow-up construction~$H$, there is some finite family~$\mathcal{F}$ of~$k$-graphs for which~$H$ is an extremal example. This then allowed him to show that~$\Pi^{(k)}_{\text{fin}}$ contains irrational numbers and that~$\Pi^{(k)}_{\infty}$ has the cardinality of the continuum. In particular, his results imply that~$\Pi^{(k)}_{\text{fin}}$ has accumulation points in $[0,1)$ for all $k \ge 3$, though we remark that the finite families given by his construction are huge. 
        In this note, we show that the same is true of~$\Pi^{(k)}$, that is, that the set of Tur\'an densities of single $k$-uniform hypergraphs has at least one accumulation point in $[0,1)$ for all $k\ge 3$. 
        



        \begin{theorem}\label{thm:mainshort}
            For every integer~$k\geq3$, the set~$\Pi^{(k)}$ has an accumulation point in~$[0,1)$.
            Moreover,~$1/2$ is an accumulation point for~$\Pi^{(3)}$.
        \end{theorem}

        This is a consequence of the following more general result.
 
        \begin{theorem}\label{thm:main}
            For every integer~$k\geq3$, there is some~$\alpha^{(k)}\in[0,1)$ such that all of the following hold:
            \begin{enumerate}
                \item There is a sequence of~$k$-graphs~$\{F_n\}_{n\in\mathds{N}}$ such that $\lim_{n\to\infty}\pi(F_n)=\alpha^{(k)}$ and~$\pi(F_n)<\alpha^{(k)}$ for all~$n\in\mathds{N}$.
                
                \item For every~$\varepsilon>0$, there is a~$k$-graph~$G_{\varepsilon}$ with~$\alpha^{(k)}\leq\pi(G_{\varepsilon})\leq \alpha^{(k)}+\varepsilon$.
                In particular,~$\pi(\{G_{\frac{1}{n}}\}_{n\in\mathds{N}})=\alpha^{(k)}$.

                \item $\alpha^{(k)}\leq\frac{k-2}{k-1}$ and~$\alpha^{(3)}=1/2$.
            \end{enumerate}
        \end{theorem}
        Perhaps surprisingly, the proofs of the first two points are abstract in the sense that they work without pinpointing~$\alpha^{(k)}$. 
        Regarding these values, it would be interesting to determine~$\alpha^{(k)}$ for $k \ge 4$ or to show that~$\alpha^{(k)}$ is itself in~$\Pi^{(k)}$ for $k \ge 3$. In particular, highlighting the depth of our ignorance about Tur\'an densities, it is not known if $1/2 \in \Pi^{(3)}$.

    \section{Preliminaries}
        Given an integer~$t$ and a~$k$-graph~$F$, let~$B(F,t)$ be the~$t$-blow-up of~$F$, the~$k$-graph obtained from~$F$ by replacing every vertex by~$t$ copies of itself.
        The following phenomenon, which we make extensive use of, is well-known.
        
        \begin{theorem}[Supersaturation]\label{thm:supersaturation}
            \begin{enumerate}
                \item\label{it:supsat:original} For every~$k$-graph~$F$ and~$\delta>0$, there are~$\varepsilon>0$ and~$n_0$ such that every~$k$-graph on~$n\geq n_0$ vertices with at least~$(\pi(F)+\delta)\binom{n}{k}$ edges contains at least~$\varepsilon n^{\vert V(F)\vert}$ copies of~$F$.
                \item\label{it:supsat:blowup} For every integer~$t$ and~$k$-graph~$F$, $\pi(B(F,t))=\pi(F)$.
                \item\label{it:supsat:hom} Let~$F$ be a~$k$-graph and let~$\mathcal{F}$ be the (finite) family of~$k$-graphs~$F'$ whose vertex set is a subset of~$V(F)$ and for which there exists a homomorphism~$\varphi:F\to F'$.
                Then~$\pi(\mathcal{F})=\pi(F)$.
                \item\label{it:supsat:onevtx} For every~$k$-graph~$F$ and~$\delta>0$, there are~$\varepsilon>0$ and~$n_0$ such that, for all~$v\in V(F)$, every~$k$-graph on~$n\geq n_0$ vertices with at least~$(\pi(F)+\delta)\binom{n}{k}$ edges contains the~$k$-graph obtained from~$F$ by replacing~$v$ by~$\varepsilon n$ copies of~$v$.
            \end{enumerate}
        \end{theorem}
    
    \section{Proof of Theorem~\ref{thm:main}}
        For~$k,\ell\in\mathds{N}$, we define the~$k$-uniform \emph{ladder of length~$\ell$} to be the~$k$-graph~$L^{(k)}_{\ell}$ with vertex set $$V(L^{(k)}_{\ell})=\{v_{ij}:i\in[\ell],j\in[k-1]\}\dcup\{t\}$$ and edge set $$E(L^{(k)}_{\ell})=\{v_{i1}\dots v_{ik-1}v_{i+1j}:i\in[\ell-1],j\in[k-1]\}\cup\{v_{\ell1}\dots v_{\ell k-1}t\}\,.$$
        For~$m\in\mathds{N}$, we further define the~$k$-graph~$L_{\ell}^{(k)}(m)$ to be the~$k$-graph with vertex set $$V(L_{\ell}^{(k)}(m))=\{v_{ij}:i\in[\ell],j\in[k-1]\}\dcup T\,,$$ where~$T$ is some set of size~$m$, and edge set $$E(L_{\ell}^{(k)}(m))=\{v_{i1}\dots v_{ik-1}v_{i+1j}:i\in[\ell-1],j\in[k-1]\}\cup\{v_{\ell1}\dots v_{\ell k-1}t:t\in T\}\,.$$
        Lastly, following~\cite{PS:23}, for an integer~$\ell\geq 2$, we define the $k$-uniform \emph{zycle of length~$\ell$} to be the~$k$-graph~$Z^{(k)}_{\ell}$ with vertex set~$V(Z^{(k)}_{\ell})=\{v_{ij}:i\in\mathds{Z}/\ell\mathds{Z},j\in[k-1]\}$ and edge set
        $$E(Z_{\ell})=\{v_{i1}\dots v_{ik-1}v_{i+1j}:i\in\mathds{Z}/\ell\mathds{Z},j\in[k-1]\}\,.$$

        Let~$k$ be an integer with~$k\geq 3$.
        We first show that there is some~$\alpha^{(k)}\in(0,1)$ such that~$\lim_{\ell\to\infty}\pi(L^{(k)}_{\ell})=\alpha^{(k)}$ while~$\pi(L^{(k)}_{\ell})<\alpha^{(k)}$ for all~$\ell\in\mathds{N}$.
        Then we show that for every~$\varepsilon>0$ there is some~$M\in\mathds{N}$ such that~$\pi(Z^{(k)}_M)\leq\alpha^{(k)}+\varepsilon$.
        Lastly, we will argue that~$\alpha^{(k)}\leq\frac{k-2}{k-1}$ with equality for~$k=3$.
        
        \subsection{Part I}
        Note that for every~$\ell\in\mathds{N}$,~$L^{(k)}_{\ell}$ is contained in a sufficiently large blow-up of~$K_{2(k-1)}^{(k)}$.
        Hence, by supersaturation (Theorem~\ref{thm:supersaturation}~\eqref{it:supsat:blowup}), we have that for all~$\ell\in\mathds{N}$,~$\pi(L^{(k)}_{\ell})\leq\pi(K_{2(k-1)}^{(k)})<1$.
        It is thus sufficient to show that~$\pi(L^{(k)}_{\ell})<\pi(L^{(k)}_{\ell+1})$ for all~$\ell\in\mathds{N}$.
        We do this by induction on~$\ell$.
        For~$\ell=1$, this follows from the result of Erd\H{o}s~\cite{E:64} mentioned in the introduction, since~$L^{(k)}_1$ is~$k$-partite, but~$L^{(k)}_2$ is not.
        Now assume that~$\ell>1$ and we know that~$\pi(L^{(k)}_i)<\pi(L^{(k)}_{i+1})$ holds for all~$i\in[\ell-1]$.
        Denote by~$\mathcal{L}^{(k)}_{\ell}$ the (finite) family of~$k$-graphs~$F$ whose vertex set is a subset of~$V(L^{(k)}_{\ell})$ and for which there exists a homomorphism~$\varphi:L^{(k)}_{\ell}\to F$.
        By supersaturation (Theorem~\ref{thm:supersaturation}~\eqref{it:supsat:hom}), we know that~$\pi(\mathcal{L}^{(k)}_{\ell})=\pi(L^{(k)}_{\ell})$.
        Therefore, it remains to show that~$\pi(\mathcal{L}^{(k)}_{\ell})<\pi(L^{(k)}_{\ell+1})$.

        Let~$\eta=\pi(L^{(k)}_{\ell})-\pi(L^{(k)}_{\ell-1})$ and note, by supersaturation (Theorem~\ref{thm:supersaturation}~\eqref{it:supsat:onevtx}), that there is some~$\varepsilon_1$ such that, for~$n$ sufficiently large, every~$k$-graph~$H$ on~$n$ vertices with at least~$(\pi(L^{(k)}_{\ell-1})+\eta/2)\binom{n}{k}$ edges contains a copy of the~$k$-graph~$L^{(k)}_{\ell-1}(\varepsilon_1n)$.
        Finally, let~$\varepsilon_2\ll\varepsilon_1,\eta$ and let~$n$ be large enough that
        \begin{align}\label{eq:constants}
            \max\Big(\Big\{\big\vert\pi(L^{(k)}_i)-\frac{\ex(n,L^{(k)}_i)}{\binom{n}{k}}\big\vert i\in\{\ell-1,\ell,\ell+1\}\Big\}\cup\Big\{\big\vert\pi(\mathcal{L}^{(k)}_{\ell})-\frac{\ex(n,\mathcal{L}^{(k)}_{\ell})}{\binom{n}{k}}\big\vert\Big\}\Big)<\varepsilon_2\,.
        \end{align}

        Now consider an extremal example~$H$ for~$\mathcal{L}^{(k)}_{\ell}$ on~$n$ vertices.
        By our choice of constants, we know that~$H$ contains a copy of~$L^{(k)}_{\ell-1}(\varepsilon_1n)$ (say with vertex set~$\{v_{ij}:i\in[\ell-1],j\in[k-1]\}\dcup T\subseteq V(H)$).
        If any~$(k-1)$-subset of~$T$ is contained in an edge of~$H$, then~$H$ would contain a (possibly) degenerate copy of~$L^{(k)}_{\ell}$, i.e., a copy of an element in~$\mathcal{L}^{(k)}_{\ell}$.
        Thus, no~$(k-1)$-subset of~$T$ is contained in an edge of~$H$.
        
        Next we add to~$H$ a complete nearly balanced $k$-partite~$k$-graph on~$T=T_1\dcup\dots\dcup T_k$ and call the resulting~$k$-graph~$H'$.
        We claim that~$H'$ is~$L^{(k)}_{\ell+1}$-free.
        Assume, for the sake of contradiction, that it contains a copy of~$L^{(k)}_{\ell+1}$ (say with vertex set~$\{u_{ij}:i\in[\ell+1],j\in[k-1]\}\dcup \{t\}\subseteq V(H)$).
        Since this copy is not contained in~$H$, one of its edges must be an edge~$x_1\dots x_k\in E(H')\setminus E(H)$, whence we also have~$x_1,\dots, x_k\in T$.
        In fact, since~$H$ is (in particular)~$L^{(k)}_{\ell}$-free, we know that, without loss of generality, there is some~$i\in[\ell]$ with~$x_1=u_{i1}\in T_1,\dots, x_{k-1}=u_{ik-1}\in T_{k-1}$.
        But then~$u_{(i+1)1},\dots,u_{(i+1)k-1}\in T_k$ and so these~$k-1$ vertices do not lie together in any edge of~$H'$, contradicting that~$\{u_{ij}:i\in[\ell+1],j\in[k-1]\}\dcup \{t\}$ is the vertex set of a copy of~$L^{(k)}_{\ell+1}$ in~$H'$.
        Hence,~$H'$ is indeed an~$L^{(k)}_{\ell+1}$-free~$k$-graph on~$n$ vertices.
        
        By~\eqref{eq:constants}, we know that~$H$ has at least~$(\pi(\mathcal{L}^{(k)}_{\ell})-\varepsilon_2)\binom{n}{k}$ edges.
        Therefore,~$H'$ has more than~$(\pi(\mathcal{L}^{(k)}_{\ell})-\varepsilon_2)\binom{n}{k}+\big(\frac{\varepsilon_1n}{k+1}\big)^k>(\pi(\mathcal{L}^{(k)}_{\ell})+\varepsilon_2)\binom{n}{k}$ edges.
        Again by~\eqref{eq:constants}, this means that~$\pi(L^{(k)}_{\ell+1})>\pi(\mathcal{L}^{(k)}_{\ell})=\pi(L^{(k)}_{\ell})$.
        We have thus proved that~$\lim_{\ell\to\infty}\pi(L^{(k)}_{\ell})=\alpha^{(k)}$  for some~$\alpha^{(k)}\in(0,1)$ with~$\pi(L^{(k)}_{\ell})<\alpha^{(k)}$ for all~$\ell\in\mathds{N}$.

        \subsection{Part II}
        Let~$\varepsilon>0$ and pick~$t,n\in\mathds{N}$ such that~$\varepsilon,k^{-1}\gg t^{-1}\gg n^{-1}$ and, for simplicity, assume that~$t\mid n$.
        Now let~$H$ be a~$k$-graph with vertex set~$[t]$ and~$e(H)\geq(\alpha^{(k)}+\varepsilon)\binom{t}{k}$.
        We will show that there is a homomorphism from~$Z^{(k)}_{\binom{t}{k-1}!}$ into~$H$.
        Let~$H_*=B(H,n/t)$ be the~$k$-graph obtained from~$H$ by replacing every vertex~$i$ of~$H$ by~$n/t$ copies of itself, the set of which we call~$V_i$.
        For~$v\in V(H_*)$, let~$f(v)$ denote the index of the partition class of~$H_*$ that contains~$v$, i.e., if~$v$ is one of the copies of a vertex~$i\in V(H)$, then~$f(v)=i$.
        Then~$H_*$ is a~$k$-graph on~$n$ vertices with~$e(H)\geq(\alpha+\varepsilon)\binom{t}{k}\big(\frac{n}{t}\big)^k\geq(\alpha^{(k)}+\varepsilon/2)\binom{n}{k}$.
        
        Since~$\pi(L^{(k)}_{\ell})<\alpha^{(k)}$ for every~$\ell$, we have that~$H_*$ contains a copy~$L$ of~$L^{(k)}_{\binom{t}{k-1}+1}$.
        As above, let~$V(L)=\{v_{ij}:i\in[\binom{t}{k-1}+1],j\in[k-1]\}\dcup\{t\}$.
        Note that for each~$i\in[\binom{t}{k-1}+1]$, the indices~$f(v_{ij})$ with~$j\in[k-1]$ are pairwise distinct, since~$v_{i1},\dots,v_{ik-1}$ are contained in an edge together.
        As~$H_*$ only has~$t$ distinct partition classes, we deduce from the pigeonhole principle that, for some~$i,i'\in[\binom{t}{k-1}+1]$ with~$i<i'$, we have~$\{f(v_{i1}),\dots,f(v_{ik-1})\}=\{f(v_{i'1}),\dots,f(v_{i'k-1})\}$.
        Since~$H_*$ is a blow-up of~$H$, this implies that there is a homomorphism of a zycle of length at most~$\binom{t}{k-1}$ into~$H$.
        As described in~\cite{PS:23}, ``cycling'' through any such zycle the right number of times yields a homomorphism from~$Z^{(k)}_{\binom{t}{k-1}!}$ into~$H$, whence~$\ex_{\text{hom}}(t,Z^{(k)}_{\binom{t}{k-1}!})\leq(\alpha^{(k)}+\varepsilon)\binom{t}{k}$.
        Since the sequence~$\ex_{\text{hom}}(m,Z^{(k)}_{\binom{t}{k-1}!})/\binom{m}{k}$ is non-increasing (in~$m$), this implies that~$\pi_{\text{hom}}(Z^{(k)}_{\binom{t}{k-1}!})\leq\alpha^{(k)}+\varepsilon$.
        Thus, we have
        \begin{align}
            \alpha^{(k)}\leq \pi(Z^{(k)}_{\binom{t}{k-1}!})=\pi_{\text{hom}}(Z^{(k)}_{\binom{t}{k-1}!})\leq\alpha^{(k)}+\varepsilon,
        \end{align}    
        where the first inequality holds since if~$H$ contains~$Z^{(k)}_{\binom{t}{k-1}!}$, then there exists a homomorphism from~$L^{(k)}_{\ell}$ into~$H$ for all~$\ell\in\mathds{N}$.
        In particular, the Tur\'an density of the family of all~$k$-uniform zycles is~$\alpha^{(k)}$.


        \subsection{Part III}
        In this subsection, we prove that~$\alpha^{(k)}\leq\frac{k-2}{k-1}$ with equality for~$k=3$.
        In~\cite{DJ:14}, DeBiasio and Jiang gave the following construction showing that~$\pi(Z^{(3)}_{\ell})\geq 1/2$ for all~$\ell\geq 2$.
        Let~$A=\{a_1,\dots,a_{\lfloor\frac{n}{2}\rfloor}\}$ and~$B=\{b_1,\dots,b_{\lceil\frac{n}{2}\rceil}\}$ and consider the~$3$-graph on~$A\cup B$ whose edges are given by all triples~$a_ib_ja_k$ and~$a_ib_jb_k$, with~$i,j<k$.
        Thus, by Part II, we have~$1/2\leq \alpha^{(3)}$.
        A result by DeBiasio and Jiang combined with an argument from~\cite{PS:23} also shows that for every~$\varepsilon>0$ there is an~$\ell$ such that~$\pi(Z^{(3)}_{\ell})\leq 1/2+\varepsilon$.
        This will also follow from Part II if we can show that~$\alpha^{(k)}\leq\frac{k-2}{k-1}$.
        By Parts I and II, it is therefore left to argue that~$\pi(L^{(k)}_{\ell})\leq\frac{k-2}{k-1}$ for every~$\ell\in\mathds{N}$. 
        Given~$\ell\in\mathds{N}$ and~$\varepsilon>0$, choose~$n\in\mathds{N}$ such that~$\varepsilon,k^{-1},\ell^{-1}\gg n^{-1}$ and let~$H$ be a~$k$-graph on~$n$ vertices with~$\delta(H)\geq (\frac{k-2}{k-1}+\varepsilon)\binom{n}{k-1}$ (a standard induction argument shows that to prove an upper bound on the Tur\'an density, we may assume such a minimum degree condition).
        Let $v_{\ell1} v_{\ell2} \dots v_{\ell k-1} t$ be an edge in $H$. Then 
        the minimum degree condition on~$H$ implies that the links of~$v_{\ell1},\dots,v_{\ell k-1}$ have at least~$\varepsilon\binom{n}{k-1}$ common edges.
        Now let~$v_{\ell-1 1},\dots,v_{\ell-1 k-1}$ be~$k-1$ vertices other than $t$ forming one of these edges.
        Using that~$\varepsilon,k^{-1}, \ell^{-1}\gg n^{-1}$, for each~$i\in[\ell-2]$ we can continue choosing vertices~$v_{i1},\dots,v_{i k-1}$ such that~$v_{i1}\dots v_{i k-1}$ is an edge in the common intersection of the links of~$v_{(i+1)1},\dots, v_{(i+1)(k-1)}$ and such that~$v_{i1},\dots, v_{i(k-1)}$ are distinct from all previously chosen vertices.
        Eventually~$v_{11},\dots,v_{1k-1},\dots,v_{\ell1},v_{\ell k-1}, t$ form the vertex set of~$L^{(k)}_{\ell}$.
        Thus,~$\pi(L^{(k)}_{\ell})\leq \frac{k-2}{k-1}$, so that, by Part I, we have~$\alpha^{(k)}\leq\frac{k-2}{k-1}$.

    \section*{Acknowledgements}
        We thank Sim\'on Piga and Marcelo Sales for interesting discussions.

\end{document}